
\newcount\secno
\newcount\prmno
\newif\ifnotfound
\newif\iffound

\def\namedef#1{\expandafter\def\csname #1\endcsname}
\def\nameuse#1{\csname #1\endcsname}

\long\def\ifundefined#1#2#3{\expandafter\ifx\csname
  #1\endcsname\relax#2\else#3\fi}
\def\hwrite#1#2{{\let\the=0\edef\next{\write#1{#2}}\next}}

\toksdef\ta=0 \toksdef\tb=2
\long\def\leftappenditem#1\to#2{\ta={\\{#1}}\tb=\expandafter{#2}%
                                \edef#2{\the\ta\the\tb}}
\long\def\rightappenditem#1\to#2{\ta={\\{#1}}\tb=\expandafter{#2}%
                                \edef#2{\the\tb\the\ta}}

\def\lop#1\to#2{\expandafter\lopoff#1\lopoff#1#2}
\long\def\lopoff\\#1#2\lopoff#3#4{\def#4{#1}\def#3{#2}}

\def\ismember#1\of#2{\foundfalse{\let\given=#1%
    \def\\##1{\def\next{##1}%
    \ifx\next\given{\global\foundtrue}\fi}#2}}

\def\section#1{\vskip1truecm
               \global\def\currenvir{section}
               \global\advance\secno by1\global\prmno=0
               {\bf \number\secno. {#1}}
               \smallskip}

\def\subsection{\global\def\currenvir{subsection}
                \global\advance\prmno by1
                \ind{(\number\secno.\number\prmno) }}
\def\subsec{\global\def\currenvir{subsection}
                \global\advance\prmno by1
                { (\number\secno.\number\prmno)\ }}

\def\proclaim#1{\global\advance\prmno by 1
                {\bf #1 \the\secno.\the\prmno$.-$ }}

\long\def\th#1 \enonce#2\endth{%
   \medbreak\proclaim{#1}{\it #2}\global\def\currenvir{th}\smallskip}

\def\bib#1{\rm #1}
\long\def\thr#1\bib#2\enonce#3\endth{%
\medbreak{\global\advance\prmno by 1\bf#1\the\secno.\the\prmno\ 
\bib{#2}$\!.-$ } {\it
#3}\global\def\currenvir{th}\smallskip}
\def\rem#1{\global\advance\prmno by 1
{\it #1} \the\secno.\the\prmno$.-$ }


\def\isinlabellist#1\of#2{\notfoundtrue%
   {\def\given{#1}%
    \def\\##1{\def\next{##1}%
    \lop\next\to\za\lop\next\to\zb%
    \ifx\za\given{\zb\global\notfoundfalse}\fi}#2}%
    \ifnotfound{\immediate\write16%
                 {Warning - [Page \the\pageno] {#1} No reference found}}%
                \fi}%
\def\ref#1{\ifx\labellist\empty{\immediate\write16
                 {Warning - No references found at all.}}
               \else{\isinlabellist{#1}\of\labellist}\fi}

\def\newlabel#1#2{\rightappenditem{\\{#1}\\{#2}}\to\labellist}
\def\labellist{}

\def\label#1{%
  \def\given{th}%
  \ifx\given\currenvir%
    {\hwrite\lbl{\string\newlabel{#1}{\number\secno.\number\prmno}}}\fi%
  \def\given{section}%
  \ifx\given\currenvir%
    {\hwrite\lbl{\string\newlabel{#1}{\number\secno}}}\fi%
  \def\given{subsection}%
  \ifx\given\currenvir%
    {\hwrite\lbl{\string\newlabel{#1}{\number\secno.\number\prmno}}}\fi%
  \def\given{subsubsection}%
  \ifx\given\currenvir%
  {\hwrite\lbl{\string%
    \newlabel{#1}{\number\secno.\number\subsecno.\number\subsubsecno}}}\fi
  \ignorespaces}

\newwrite\lbl

\def\openall{\openout\lbl=\jobname.lbl}

\newread\testfile
\def\lookatfile#1{\openin\testfile=\jobname.#1
    \ifeof\testfile{\immediate\openout\nameuse{#1}\jobname.#1
                    \write\nameuse{#1}{}
                    \immediate\closeout\nameuse{#1}}\fi%
    \immediate\closein\testfile}%

\def\begin{\newlabel{def*}{1.1}
\newlabel{grad}{1.2}
\newlabel{defF}{1.3}
\newlabel{Finv}{1.4}
\newlabel{F*}{1.5}
\newlabel{level}{1.6}
\newlabel{FF}{1.7}
\newlabel{newt}{2.2}
\newlabel{FC}{2.3}
\newlabel{Nk}{2.4}
\newlabel{bigrad}{2.5}
\newlabel{stabF}{3.1}
\newlabel{d}{4.1}
}
           

\magnification 1250
\pretolerance=500 \tolerance=1000  \brokenpenalty=5000
\mathcode`A="7041 \mathcode`B="7042 \mathcode`C="7043
\mathcode`D="7044 \mathcode`E="7045 \mathcode`F="7046
\mathcode`G="7047 \mathcode`H="7048 \mathcode`I="7049
\mathcode`J="704A \mathcode`K="704B \mathcode`L="704C
\mathcode`M="704D \mathcode`N="704E \mathcode`O="704F
\mathcode`P="7050 \mathcode`Q="7051 \mathcode`R="7052
\mathcode`S="7053 \mathcode`T="7054 \mathcode`U="7055
\mathcode`V="7056 \mathcode`W="7057 \mathcode`X="7058
\mathcode`Y="7059 \mathcode`Z="705A
\def\spacedmath#1{\def\packedmath##1${\bgroup 
\mathsurround =0pt##1\egroup$}\mathsurround#1
\everymath={\packedmath}\everydisplay={\mathsurround=0pt}}
\def\nospacedmath{\mathsurround=0pt
\everymath={}\everydisplay={} } \spacedmath{2pt}
\def\qfl#1{\buildrel {#1}\over {\longrightarrow}}
\def\phfl#1#2{\normalbaselines{\baselineskip=0pt
\lineskip=10truept\lineskiplimit=1truept}\nospacedmath\smash 
{\mathop{\hbox to 8truemm{\rightarrowfill}}
\limits^{\scriptstyle#1}_{\scriptstyle#2}}}
\def\hfl#1#2{\normalbaselines{\baselineskip=0truept
\lineskip=10truept\lineskiplimit=1truept}\nospacedmath
\smash{\mathop{\hbox to
12truemm{\rightarrowfill}}\limits^{\scriptstyle#1}_{\scriptstyle#2}}}
\def\diagram#1{\def\normalbaselines{\baselineskip=0truept
\lineskip=10truept\lineskiplimit=1truept}   \matrix{#1}}
\def\vfl#1#2{\llap{$\scriptstyle#1$}\left\downarrow\vbox to
6truemm{}\right.\rlap{$\scriptstyle#2$}}
\def\mono{\lhook\joinrel\mathrel{\longrightarrow}}
\def\iso{\vbox{\hbox to .8cm{\hfill{$\scriptstyle\sim$}\hfill}
\nointerlineskip\hbox to .8cm{{\hfill$\longrightarrow $\hfill}} }}

\def\sdir_#1^#2{\mathrel{\kern3pt\mathop{\oplus}\limits_{#1}^{#2}}}
\def\pprod_#1^#2{\raise2pt\hbox{$\mathrel{\scriptstyle\mathop
{\kern0pt\prod}\limits_{#1}^{#2}}$}}

\font\eightrm=cmr8         \font\eighti=cmmi8
\font\eightsy=cmsy8        \font\eightbf=cmbx8
\font\eighttt=cmtt8        \font\eightit=cmti8
\font\eightsl=cmsl8        \font\sixrm=cmr6
\font\sixi=cmmi6           \font\sixsy=cmsy6
\font\sixbf=cmbx6\catcode`\@=11
\def\eightpoint{%
  \textfont0=\eightrm \scriptfont0=\sixrm \scriptscriptfont0=\fiverm
  \def\rm{\fam\z@\eightrm}%
  \textfont1=\eighti  \scriptfont1=\sixi  \scriptscriptfont1=\fivei
  \def\oldstyle{\fam\@ne\eighti}\let\old=\oldstyle
  \textfont2=\eightsy \scriptfont2=\sixsy \scriptscriptfont2=\fivesy
  \textfont\itfam=\eightit
  \def\it{\fam\itfam\eightit}%
  \textfont\slfam=\eightsl
  \def\sl{\fam\slfam\eightsl}%
  \textfont\bffam=\eightbf \scriptfont\bffam=\sixbf
  \scriptscriptfont\bffam=\fivebf
  \def\bf{\fam\bffam\eightbf}%
  \textfont\ttfam=\eighttt
  \def\tt{\fam\ttfam\eighttt}%
  \abovedisplayskip=9pt plus 3pt minus 9pt
  \belowdisplayskip=\abovedisplayskip
  \abovedisplayshortskip=0pt plus 3pt
  \belowdisplayshortskip=3pt plus 3pt 
  \smallskipamount=2pt plus 1pt minus 1pt
  \medskipamount=4pt plus 2pt minus 1pt
  \bigskipamount=9pt plus 3pt minus 3pt
  \normalbaselineskip=9pt
  \setbox\strutbox=\hbox{\vrule height7pt depth2pt width0pt}%
  \normalbaselines\rm}\catcode`\@=12

\newcount\noteno
\noteno=0
\def\up#1{\raise 1ex\hbox{\sevenrm#1}}
\def\note#1{\global\advance\noteno by1
\footnote{\parindent0.4cm\up{\number\noteno}\
}{\vtop{\eightpoint\baselineskip12pt\hsize15.5truecm\noindent
#1}}\parindent 0cm}

\def\pc#1{\tenrm#1\sevenrm}
\def\tx{\kern-1.5pt -}
\def\cqfd{\kern 2truemm\unskip\penalty 500\vrule height 4pt depth 0pt width
4pt} 
\def\virg{\raise
.4ex\hbox{,}}
\def\decale#1{\smallbreak\hskip 28pt\llap{#1}\kern 5pt}
\def\no{n\up{o}\kern 2pt}
\def\ind{\par\hskip 1truecm\relax}
\def\indp{\par\hskip 0.5truecm\relax}

\def\rond{\kern 1pt{\scriptstyle\circ}\kern 1pt}

\def\det{\mathop{\rm det}\nolimits}

\def\ch{\mathop{\rm ch}\nolimits}
\def\td{\mathop{\rm Todd}\nolimits}

\font\gragrec=cmmib10
\def\varphig{\hbox{\gragrec \char39}}
\font\ptigrec=cmmib7
\def\varphip{\hbox{\ptigrec\char39}}

\frenchspacing
\input amssym.def
\input amssym
\vsize = 25truecm
\hsize = 16truecm
\voffset = -.5truecm
\parindent=0cm
\baselineskip15pt
\overfullrule=0pt

\begin
\centerline{\bf Algebraic cycles on Jacobian varieties}
\smallskip
\smallskip \centerline{Arnaud {\pc BEAUVILLE}} 
\vskip1.2cm

{\bf Introduction}
\smallskip
\ind Let $C$ be a compact Riemann surface of genus $g$. Its
Jacobian variety $J$ carries a number
of natural subvarieties, defined up to translation: first of all the
curve $C$ embeds  into $J$, then we can use the group
law in
$J$ to form   $W_2=C+C$, $W_3=C+C+C$, ... till
$W_{g-1}$ which is a  theta
divisor  on $J$. Then we can intersect these
subvarieties, add again, pull back or push down under
multiplication by integers, and so on. Thus we get a rather big
supply of algebraic subvarieties which live naturally in $J$.
\ind If we look at the classes  obtained  this way in
rational cohomology, the result is disappointing: we just
find the subalgebra of $H^*(J,{\bf Q})$ generated by the class
$\theta$ of the theta divisor -- in fact, the polynomials in $\theta$ are the
only algebraic cohomology classes which live on a generic
Jacobian. The situation becomes more interesting if we look at
the
${\bf Q}$\tx algebra $A(J)$ of algebraic cycles modulo algebraic
equivalence on $J$; here a result of Ceresa  [C] implies that for a
generic curve $C$, the  class of $W_{g-p}$ in $A^p(J)$ is
{\it not} proportional to $\theta^p$   for $2\le p\le g-1$. This leads naturally
to investigate the ``tautological subring" of $A(J)$, that is,
the smallest ${\bf Q}$\tx vector subspace $R$ of
$A(J)$ which contains $C$ and is stable under the natural
operations of $A(J)$: intersection and Pontryagin products (see
(\ref{def*}) below), pull back and push down under
multiplication by integers. Our main result states that this space is
not too complicated. Let $w^p\in A^p(J)$ be the class of $W_{g-p}$.
Then:\smallskip 
 {\bf  Theorem}$.-$ a) {\it $R$ is the sub-$\!{\bf Q}$\tx algebra
of $A(J)$ generated by} $w^1,\ldots,w^{g-1} $.
\ind b) {\it If $C$ admits a morphism of degree $d$ onto ${\bf
P}^1$, $R$ is generated by} $w^1,\ldots,w^{d-1}$.
\smallskip 
\ind In particular we see that $R$ is finite-dimensional, a fact
which does not seem to be  a priori obvious (the space $A(J)$ is
known to be infinite-dimensional for $C$ generic of genus $3$, see [N]). 
\ind The proof rests in an essential way on the properties of the
Fourier transform, a ${\bf Q}$\tx linear automorphism of 
$A(J)$ with remarkable properties. We recall  these 
properties in
\S 1; in \S 2 we look at the case of
Jacobian varieties, computing in particular the Fourier transform
of the class of $C$ in $A(J)$. This is the main ingredient  in the
proof of part a) of the theorem, which we give in \S 3. Part b)
turns out to be an easy consequence of a result of Colombo and
van Geemen [CG]; this is explained in
\S 4, together with a few examples.

\section{Algebraic cycles on abelian varieties}
\subsection 
Let $X$ be an abelian variety over ${\bf C}$. We will denote by
$p$ and $q$ the two projections of $X\times X$ onto $X$, and  by
$m:X\times X\rightarrow X$ the addition map.
\ind Let $A(X)$ be the group of
algebraic cycles on $X$ modulo algebraic equivalence, {\it
tensored with} ${\bf Q}$. It is a ${\bf Q}$\tx vector space,
graded by the codimension of the cycle classes. It carries two
natural multiplication laws $A(J)\otimes_{\bf Q}A(J)\rightarrow
A(J)$, which are associative and commutative:  the
intersection product, which is homogeneous with respect to the
graduation, and the Pontryagin product, defined by
$$x*y:=m_*(p^*x\cdot q^*y)\ ,$$
which is homogeneous of degree $-g$. If $Y$ and $Z$ are
subvarieties of $X$, the cycle class $[Y]*[Z]$ is equal to $(\deg
\mu)\,[Y+Z]$ if the addition map $\mu:Y\times Z\rightarrow Y+Z$ is
generically finite, and is zero otherwise.\label{def*}\smallskip 
\subsection For $k\in {\bf Z}$, we will still denote by $k$ the
endomorphism $x\mapsto kx$ of $X$. According to [B2], there is
a second graduation on $A(X)$, leading to a bigraduation
$$A(X)= \sdir_{s,p}^{}A^p(X)_{(s)}$$
such that
$$k^*x=k^{2p-s}x\quad,\quad k_*x=k^{2g-2p+s}x\quad{\rm
for}\ x\in A^p(X)_{(s)}\ .$$
Both products are homogeneous with respect to the second graduation. We
have $A^p(X)_{(s)}=0$ for $s<p-g$ or
$s\ge g$ (use [B2], Prop. 4). It is conjectured that negative degrees
actually do not occur; this will not concern us here, as  we
will  only consider cycles in $A(X)_{(s)}$ for $s\ge 0$. \label{grad}
\smallskip 
\subsection\label{defF}
A crucial tool in what follows will be the Fourier
transform for algebraic cycles, defined in [B1]. Let us recall
briefly the results  we will need  -- the proofs can
be found in [B1] and [B2].  We will concentrate on
the case of a {\it principally polarized} abelian variety
$(X,\theta)$, identifying $X$ with its dual abelian variety.  
\ind  Let $\ell :=p^*\theta+q^*\theta-m^*\theta\in A^1(X\times
X)$; it is the class of the {\it  Poincar\'e line bundle} ${\cal L}$
on $X\times X$. 
  The Fourier transform
${\cal F}:A(X)\rightarrow A(X)$ is defined by ${\cal
F}x=q_*(p^*x\cdot e^\ell )$. It satisfies the following
properties:
\indp\subsec \label{Finv} ${\cal F}\rond {\cal
F}=(-1)^g(-1)^*$
\indp\subsec \label {F*}  ${\cal F}(x*y)={\cal F}x\cdot {\cal F}y$ 
and  ${\cal F}(x\cdot y)=(-1)^g\,{\cal F}x*{\cal F}y$
\indp\subsec \label{level} ${\cal
F}A^p(X)_{(s)}=A^{g-p+s}(X)_{(s)}$
\indp\subsec \label{FF} Let $x\in A(X)$; put $\bar
x=(-1)^*x$. Then ${\cal F}x=e^\theta\,\bigl ((\bar x
e^\theta)*e^{-\theta}\bigr)$.
\ind Let us prove (\ref{FF}), which is not explicitly stated
in [B1] or [B2].  Replacing $\ell $ by its definition, we
get
${\cal F}x=e^\theta\, q_*(p^*(xe^\theta)\cdot e^{-m^*\theta})$.
Let  $\omega$ be the automorphism of $A\times A$ defined
by $\omega(a,b)=(-a,a+b)$. We have $p\rond\omega=-p$,
$q\rond\omega=m$, $m\rond\omega=q$. Hence
$${\cal F}x=e^\theta\,
q_*\omega_*\omega^*(p^*(xe^\theta)\cdot
e^{-m^*\theta})=e^\theta\, m_*\bigl(p^*(\bar x e^\theta)\cdot
q^*e^{-\theta}\bigr) = e^\theta\,\bigl ((\bar x
e^\theta)*e^{-\theta}\bigr)\ .\cqfd
$$
\section{The Fourier transform on a Jacobian}
\subsection  From now on we take for our abelian variety the
Jacobian 
$(J,\theta)$ of a smooth projective curve $C$ of genus
$g$. We choose a base point $o\in C$, which allows us to define
an embedding $\varphi:C\mono J$ by $\varphi(p)={\cal
O}_C(p-o)$. Since we are working modulo algebraic equivalence,
all our constructions will be independent of the choice of the base
point.
\ind We will denote simply by $C$ the class of $\varphi(C)$ in
$A^{g-1}(J)$. For $0\le d\le g$, we put $w^{g-d}:={1\over d!}C^{*d}\in
A^{g-d}(J)$; it is the class of the subvariety $W_d$ of $J$ parameterizing line
bundles of the form
${\cal O}_C^{}
(E_d-do)$, where
$E_d$ is an effective divisor of degree $d$. We have
$w^1=\theta$ by the Riemann theorem, $w^{g-1}=C$, and $w^g$ is the class
of a point. We define the {\it Newton polynomials} in the classes $w^i$ by 
$$N^k(w)={1\over k!}\sum_{i=1}^g \lambda_i^k$$
in the ring obtained by adjoining to $A(J)$ the roots
$\lambda_1,\ldots,\lambda_g$ of the equation
$\lambda^n-\lambda^{n-1}\,w^1+\ldots+(-1)^gw^g=0$. We have
$N^k(w)\in A^k(J)$; for instance:
$$N^1(w)=\theta\quad ,\quad 
N^2(w)={1\over 2}\theta^2-w^2\quad ,\quad N^3(w)={1\over
6}\theta^3-{1\over 2}\theta\cdot w^2 -{1\over 2}w^3\quad
,\quad  ...
$$
\subsection The class $N^k(w)$ is a polynomial in $w^1,\ldots,
w^k$; conversely,  
$w^k$ is a polynomial in
$N^1(w),\ldots ,N^k(w)$.\label{newt}
\th Proposition
\enonce We have $\  -{\cal
F}C=N^1(w)+N^2(w)+\ldots+ N^{g-1}(w)$.
\endth\label{FC}
{\it Proof} : We use the notation of (\ref{defF}), and denote
moreover by 
$\bar p,\bar q$  the  projections of
$C\times J$ onto $C$ and $J$. 
Consider the cartesian diagram
$$\diagram{C\times J&\hfl{\Phi}{}& J\times J\cr
\vfl{\bar p}{}&&\vfl{}{ p}\cr
C & \hfl{}{\varphi} & J
}$$
with $\Phi=(\varphi,1_J)$. Put $\bar\ell :=\Phi^*\ell $. We have
$ p^*C \cdot e^{\ell }=\Phi_*1\cdot e^\ell
=\Phi_*e^{\bar\ell }$, and therefore 
$${\cal F}C=\bar q_*e^{\bar\ell }\ .$$
\ind The line bundle $\bar{\cal L}:=\Phi^*{\cal L}$ is the
Poincar\'e line bundle on $C\times J$: that is, we have $\bar {\cal
L}_{C\times \{\alpha\}}=\alpha$ for all $\alpha\in J$, and $\bar
{\cal L}_{\{o\}\times J}={\cal O}_J$. We will now work
exclusively on $C\times J$, and suppress the bar above the
letters $p$, $q$, ${\cal L}$ and $\ell $. We  apply the
Grothendieck-Riemann-Roch theorem to
$ q$ and
${\cal L}$. Since we are working modulo algebraic
equivalence, the Todd class of
$C$ is simply $1+(1-g)o$. Let $i_o:J\mono C\times J$ be the map
$\alpha\mapsto (o,\alpha)$; we have
$$  q_*( p^*o\cdot
e^{\ell })=q_*i_{o*}i_o^*e^\ell =i_o^*e^\ell =1$$
since $i_o^*{\cal L}$ is trivial. Thus
$$\ch q_!{\cal L}= q_*(p^*\td(C)\cdot \ch {\cal L})=q_*e^\ell 
-(g-1)\ .$$
\ind  The Chern classes of 
$q_!{\cal L}$ are computed in [M]: we have
$$c(-q_!{\cal L})=1+w^1+\ldots+w^g\ .$$
Putting things together we obtain
$${\cal F}C=q_*e^\ell  =g-1-\ch (-q_!{\cal
L})=-(N^1(w)+N^2(w)+\ldots+N^g(w))\ .\cqfd$$
\ind Let $\displaystyle C=\sum_{s=0}^{g-1}C_{(s)}$ be the
decomposition of $C$ in $\sdir_s^{} A^{g-1}(J)_{(s)}$.  From the
Proposition,   (\ref{level}) and (\ref{Finv})  we obtain:
\th Corollary
\enonce We have $N^k(w)=-{\cal F}C_{(k-1)}\in A^k(J)_{(k-1)}$
and
${\cal F}(N^k(w))=$ $(-1)^{g+k}C_{(k-1)}\,.\cqfd$ 
\endth\label{Nk}
\th Corollary
\enonce The ${\bf Q}$\tx subalgebra $R$ of $A(J)$ generated by 
$w^1,\ldots,w^{g-1}$ is  bigraded. In particular, it is stable under
the operations $k^*$ and $k_*$ for each $k\in {\bf Z}$.
\endth\label{bigrad}
\ind Indeed $R$ is also generated by the elements $N^1(w),\ldots
,N^{g-1}(w)$ (\ref{newt}), which are homogeneous for both
graduations.\cqfd
\section{Proof of the main result}
\ind In order to prove part a) of the theorem, it remains to prove
that the \hbox{${\bf Q}$\tx sub\-algebra} $R$ of $A(J)$ generated
by 
$w^1,\ldots,w^{g-1}$ is stable under the Pontryagin product. In
view of (\ref{F*}) it suffices to prove:
\th Proposition
\enonce $R$  is stable under ${\cal F}$.
\endth\label{stabF}
\smallskip 
{\it Proof} : Let ${\cal F}R$ denote the image of $R$
under the Fourier transform; it is a vector space over ${\bf Q}$, 
stable under the Pontryagin product (\ref{F*}).  We will prove that
${\cal F}R$ is stable under ${\cal F}$, that is,   ${\cal F}{\cal
F}R\i {\cal F}R$; since ${\cal F}{\cal F}R=R$ (\ref{Finv}) this
implies 
$R\i{\cal F}R$, then ${\cal F}R\i R$ by applying ${\cal F}$ again.
 \ind We observe that is enough to prove that ${\cal F}R$ {\it is
stable under multiplication by} $\theta$. Indeed it is then stable
under multiplication by $e^{\theta}$, and finally under ${\cal
F}$ in view of the formula ${\cal
F}x=e^\theta\,\bigl ((\bar x e^\theta)*e^{-\theta}\bigr)$
(\ref{FF}).
\ind Since the ${\bf Q}$\tx algebra $R$ is generated by 
 the classes
$N^p(w)$, ${\cal F}R$ is spanned as a ${\bf Q}$\tx vector space
by the elements
$${\cal F}(N^{p_1}(w)\ldots N^{p_r}(w))=\pm\, C_{(p_1-1)}*
\ldots *C_{(p_r-1)}$$
(we are using (\ref{F*}) and Corollary \ref{Nk}).
Actually:
\th Lemma 
\enonce ${\cal F}R$ is spanned by the classes $(k_{1*}C)*\ldots
*(k_{r*}C)$, for all sequences $(k_1,\ldots ,k_r)$ of positive
integers.
\endth
{\it Proof} : For
$k\in {\bf Z}$ we have $\displaystyle
k_*C=\sum_{s=0}^{g-1}k^{2+s}C_{(s)}$ (\ref{grad}); therefore
$$(k_{1\,*}C)*\ldots*(k_{r\,*}C)=(k_1\ldots
k_r)^2\sum_{s_1,\ldots,s_r}k_1^{s_1}\ldots k_r^{s_r}\
C_{(s_1)*}\ldots *C_{(s_r)}
$$
where ${\bf s}=(s_1,\ldots,s_r)$ runs in $[0,g-1]^r$; this shows in
particular that  $(k_{1\,*}C)*\ldots$ $*(k_{r\,*}C)$ belongs to
${\cal F}R$. We claim that we can choose $g^r$ $r$\tx uples 
${\bf k}=(k_1,\ldots,k_r)$ so that the matrix $(a_{{\bf
k},{\bf s}})$ with entries $a_{{\bf k},{\bf s}}=
(k_1^{s_1}\ldots k_r^{s_r})$ is
invertible: if we take for instance the sequence of $r$\tx
uples ${\bf k}_\ell =(\ell ,\ell ^{g},\ldots, \ell
^{g^{r-1}})$, for
$1\le\ell \le g^r$, we get for $\det(a_{{\bf
k},{\bf s}})$ a non-zero Vandermonde determinant. Thus 
each element $C_{(s_1)*}\ldots
*C_{(s_r)}$ is a  ${\bf Q}$\tx linear combination of 
classes of the form $(k_{1\,*}C)*\ldots
*(k_{r\,*}C)$, which proves the lemma.\cqfd\smallskip

\subsection Thus it suffices  to prove that 
each  product
$\theta\cdot \bigl((k_{1*}C)*\ldots
*(k_{r*}C)\bigr)$ belongs to ${\cal F}R$.
We observe that $(k_{1*}C)*\ldots
*(k_{r*}C)$ is a multiple of the image of the composite map
$$u:C^r\qfl{\varphip} J^r\qfl{{\bf k}}J^r\qfl{m}J $$ 
where ${\bf k}=(k_1,\ldots,k_r)$, $\varphig=(\varphi,\ldots
,\varphi)$ and $m$ is the addition morphism.  Thus the class 
$\theta\cdot \bigl((k_{1*}C)*\ldots
*(k_{r*}C)\bigr)$ is proportional to
$u_*u^*\theta$.
\ind Let $p_i: J^r\rightarrow J$ (resp. $p_{ij}:J^r\rightarrow
J^2$) denote the projection onto the $i$\tx th factor (resp. the
$i$\tx th and $j$\tx th factors).  In $A^1(J^r)$, we have
$$m^*\theta=\sum_i p_i^*\theta-\sum_{i<j}p_{ij}^*\ell \
;$$indeed for $r=2$ this is the definition of
$\ell $, and the general case follows from the theorem of the cube.
We have also
$k_i^*\theta=k_i^2\,\theta$ and
$(k_i,k_j)^*\ell =k_ik_j\,\ell $. Thus
$${\bf
k}^*m^*\theta=\sum_ik_i^2\,p_i^*\theta-\sum_{i<j}k_ik_j\,
p_{ij}^*\ell \ ;$$denoting  by $q_i$, $q_{ij}$ the
projections of $C^r$ onto $C$ and $C^2$, we find
$$u^*\theta=\sum_ik_i^2\,q_i^*\varphi^*\theta-\sum_{i<j}k_ik_j\,
q_{ij}^*(\varphi,\varphi)^*\ell \ .$$
Let $\Delta$ be the diagonal in $C^2$. The theorem of the
square gives  $$(\varphi,\varphi)^*{\cal L}={\cal
O}_{C^2}(\Delta-C\times o\,-\,o\times C)\ .$$ Therefore
$u^*\theta$ is algebraically equivalent to a linear combination of
divisors of the form $q_i^*o$ and
$q_{ij}^*\Delta$. Under $u_*$ each of these divisors is mapped
to a multiple of the cycle
$(l_{1*}C)*\ldots (l_{r-1*}C)$, where the sequence
$(l_{1}\ldots l_{r-1})$ is $(k_1,\ldots,\widehat{k_i},\ldots
,k_r)$ in the first case and
$(k_1,\ldots,\widehat{k_i},\ldots,\widehat{k_j},\ldots  ,k_r,
k_i+k_j)$ in the second one (as usual the symbol
$\widehat{k_i}$ means that $k_i$ is omitted). This proves our
claim, and therefore the Proposition.\cqfd

\section{$\hbox{\gragrec\char100}$\tx gonal curves}
\th Proposition
\enonce Assume that the curve $C$ is $d$\tx gonal, that is, admits
a degree $d$ morphism onto ${\bf P}^1$. We have $N^k(w)=0$
for $k\ge d$, and the ${\bf Q}$\tx algebra $R$ is generated by
$w^1,\ldots, w^{d-1}$.
\endth\label{d}
{\it Proof} : By now this is an immediate consequence of a result
of Colombo and van Geemen, which says that for a $d$\tx
gonal curve  $C_{(s)}=0$ for
$s\ge d-1$ ([CG], Prop. 3.6)\note{Our class $C_{(s)}$ is denoted
 $\pi _{2g-2-s}C$ in [CG]}.  This implies
$N^k(w)=0$ for
$k\ge d$ (Prop. \ref{Nk}), so that $R$ is a polynomial ring in
$N^1(w),\ldots, N^{d-1}(w)$, hence in $w^1,\ldots, w^{d-1}$
(\ref{newt}).\cqfd\medskip  
\ind The case $d=2$ of the Proposition had already
been observed by Collino [Co]: 
\th Corollary
\enonce If $C$ is hyperelliptic, $R={\bf
Q}[\theta]/(\theta^{g+1})$.
\endth
\th Corollary
\enonce If $C$ is trigonal, $R$ is generated by $\theta$ and
the class $\eta=N^2(w)$ in $A^2(J)$. There exists an integer $k\le
{g\over 3}$ such that $$R={\bf
Q}[\theta,\eta]/(\theta^{g+1},\theta^{g-2}\eta,\ldots,
\theta^{g+1-3k}\eta^k,\eta^{k+1})\ .$$
\endth
{\it Proof} : By Proposition (\ref{d}) $R$ is generated by $\theta$
and
$\eta$. For $p,s\in {\bf N}$, the class
$\theta^{p-2s}\eta^s$ is the only  monomial in $\theta,\eta$
which belongs to  $A^{p}(J)_{(s)}$; therefore it spans the
${\bf Q}$\tx vector space $R^{p}_{(s)}$ (in particular, this space
is zero for 
$p<2s$). This implies  that the relations between
$\theta$ and
$\eta$ are monomial, that is, of the form $\theta^r\eta^s=0$ for
some pairs $(r, s)\in{\bf N}^2$.
\ind Similarly, as a ${\bf Q}$\tx algebra for the Pontryagin
product, $R$ is generated by $C_{(0)}$ and $C_{(1)}$. The ${\bf
Q}$\tx vector space  $R^{p}_{(s)}$ is 
 spanned by $C_{(0)}^{*(g-p-s)}*C_{(1)}^{*s}$, hence is zero for
$p+s>g$. In particular we see that $\theta^r\eta^s=0$ as soon as
$r+3s>g$. 
\ind Let $k$ be the smallest integer such that $\eta^k\not=0$,
$\eta^{k+1}=0$. By what we have just seen the first relation
implies $3k\le g$. Suppose we have $\theta^r\eta^s=0$ for some
integers $r,s$ with $r+3s\le g$ and $s\le k$. Then we have
$R^{r+2s}_{(s)}=0$ and $C_{(0)}^{*(g-r-3s)}*C_{(1)}^{*s}=0$.
Taking $*$\tx product with $C_{(0)}^{*r}$ we arrive at
 $C_{(0)}^{*(g-3s)}*C_{(1)}^{*s}=0$, which implies $\eta^s=0$,
contradicting the definition of $k$.\cqfd\smallskip 
\ind In the general case, since any curve of genus $g$ has a
$g^1_d$ with $d\le {g+3\over 2}$ ([ACGH], Ch. V, thm. 1.1), we get:
\th Corollary
\enonce Put $ d:=[{g+1\over 2}]$. The ${\bf Q}$\tx
algebra $R$ is generated by $w^1,\ldots, w^{d}$.\cqfd
\endth\smallskip 
\subsection On the opposite side, we may ask how many of the
classes $w^i$ are really needed to generate $R$. Here we know
nothing except Ceresa's result [C]:
 {\it for $C$ generic of genus $\ge 3$ the class $N^2(w)$ is
non-zero}. Ceresa uses the intermediate Jacobian of $J$,
thus his method can only detect non-zero classes in $A(J)_{(1)}$.
It would be quite interesting to decide for instance whether 
$N^3(w)$ is nonzero for a generic curve of genus $\ge 5$.

\vskip2cm
\centerline{ REFERENCES} \vglue15pt\baselineskip12.8pt
\def\num#1{\smallskip\item{\hbox to\parindent {\enskip
[#1]\hfill}}}
\parindent=1cm 
\num{ACGH}\hskip9pt  E. {\pc ARBARELLO}, M. {\pc CORNALBA}, P. {\pc
GRIFFITHS}, J. {\pc HARRIS}: {\sl Geometry of alge\-braic curves} I. Grund. der
math. Wiss. {\bf 267}, Springer-Verlag, New York-Berlin-Heidelberg-Tokyo
(1985). 

\num{B1} A. {\pc BEAUVILLE}: {\sl   Quelques remarques sur la
transformation de Fourier dans l'anneau de Chow d'une
vari\'et\'e ab\'elienne}.  Algebraic Geometry (Tokyo/Kyoto
1982), LN {\bf 1016},  238--260; Springer-Verlag (1983). 
\num{B2} A. {\pc BEAUVILLE}: {\sl 	Sur l'anneau de Chow d'une
vari\'et\'e ab\'elienne}. Math. Ann. {\bf 273}  (1986),
647--651.

\num{C} G. {\pc CERESA}: {\sl $C$ is not algebraically equivalent 
to $C\sp{-}$ in its Jacobian}. Ann. of Math. (2) {\bf 117} (1983),
285--291.
\num{CG} E. {\pc COLOMBO}, B. {\sevenrm VAN} {\pc GEEMEN}:
{\sl Note on curves in a Jacobian}. Compositio Math. {\bf 88}
(1993),  333--353.
\num{Co} A. {\pc COLLINO}: {\sl Poincar\'e's formulas and
hyperelliptic curves}. Atti Accad. Sci. Torino {\bf 109} (1975),
89--101.
\num{M} A. {\pc MATTUCK}: {\sl
Symmetric products and Jacobians}.
Amer. J. Math. {\bf 83} (1961), 189--206
\num{N} M. {\pc NORI}: {\sl Cycles on the generic abelian
threefold}. Proc. Indian Acad. Sci. Math. Sci. {\bf 99} (1989), 
191--196.

\vskip1cm
\def\pc#1{\eightrm#1\sixrm}
\hfill\vtop{\eightrm\hbox to 5cm{\hfill Arnaud  {\pc
BEAUVILLE}\hfill}
 \hbox to 5cm{\hfill Laboratoire J.-A. Dieudonn\'e\hfill}
 \hbox to 5cm{\sixrm\hfill UMR 6621 du CNRS\hfill}
\hbox to 5cm{\hfill {\pc UNIVERSIT\'E DE}  {\pc NICE}\hfill}
\hbox to 5cm{\hfill  Parc Valrose\hfill}
\hbox to 5cm{\hfill F-06108 {\pc NICE} Cedex 02\hfill}}
\end